\documentclass[12pt]{article}
\usepackage{indentfirst} 
\usepackage{amsfonts,amssymb,amsmath,amsthm}
\usepackage{bm}          
\usepackage{cite}
\usepackage{url}
\usepackage{enumitem}
\usepackage[usenames,dvipsnames]{xcolor}
\usepackage{todonotes}
\usepackage{multicol}
\usepackage{MnSymbol}
\usepackage{hyperref}
\usepackage{textcomp}
\usepackage[titletoc,title]{appendix}
\usepackage{smartref}  \addtoreflist{section}   

\theoremstyle{definition}
\newtheorem*{propertiess}{Properties}

\newcommand{\refconchap}[1]{\hyperref[#1]{\textcolor{\corlink}{\sectionref{#1}.}\ref*{#1}}}   

\newcommand{\margenizquierda}{2.8cm}   
\newcommand{\margenderecha}{2.8cm}     
\newcommand{\margensuperior}{3cm}      
\newcommand{\margeninferior}{3.5cm}    
\usepackage{geometry} \allowdisplaybreaks

\makeatletter
\newcommand{\raisemath}[1]{\mathpalette{\raisem@th{#1}}}
\newcommand{\raisem@th}[3]{\raisebox{#1}{$#2#3$}}
\makeatother

\def\QED{{\boldmath$\rule{0.5em}{0.5em}$}}
\def\markatright#1{\leavevmode\unskip\nobreak\quad\hspace*{\fill}{#1}}
\def\qed{\markatright{\QED}}

\newtheorem{theorem}{Theorem}[section]

\newtheorem{corollary}[theorem]{Corollary}

\newcommand{\N}{\mathbb{N}}     
\newcommand{\Z}{\mathbb{Z}}     
\newcommand{\R}{\mathbb{R}}     
\newcommand{\C}{\mathbb{C}}     

\newcommand{\corurl}{blue}
\newcommand{\corcite}{ForestGreen}
\newcommand{\corlink}{blue}

\hypersetup{linktocpage,colorlinks,urlcolor=\corurl,citecolor=\corcite,%
            linkcolor=\corlink, pdftitle={Appell polynomials},
	    pdfauthor={J.F. Barbero, J. Salas, E.J.S. Villasen\~or}}

\newcommand{\arxiv}[1]{\href{http://arxiv.org/abs/#1}{\texttt{arXiv:#1}}}

\usepackage{titling}
\usepackage{authblk}
\title{On the asymptotics of the rescaled Appell polynomials}
\author[1,3]{J. Fernando Barbero G.}
\author[2,3]{Jes\'us Salas\,}
\author[2,3]{Eduardo J.S. Villase\~nor\,}
\affil[1]{Instituto de Estructura de la Materia, CSIC. Serrano 123,
          28006 Madrid, Spain}
\affil[2]{Universidad Carlos III de Madrid. Avda.\  de la Universidad
          30, 28911 Legan\'es, Spain}
\affil[3]{Grupo de Teor\'{\i}as de Campos y F\'{\i}sica Estad\'{\i}stica.
          Instituto Gregorio Mill\'an (UC3M). Unidad Asociada al Instituto
          de Estructura de la Materia, CSIC, Spain}
\date{}                     
\date{April 30, 2019}

\numberwithin{equation}{section}
\begin{document}

\maketitle

\begin{abstract}
\noindent We introduce a new representation for the rescaled Appell
polynomials and use it to obtain asymptotic expansions to 
\textit{arbitrary} order.
This representation consists of a finite sum and an integral over a universal
contour (i.e. independent of the particular polynomials 
considered within the Appell
family). We illustrate our method by studying the zero attractors for
rescaled Appell polynomials. We also discuss the asymptotics to
arbitrary order of the rescaled Bernoulli polynomials.
\end{abstract}

\bigskip

\noindent
{\bf Keywords}: Appell polynomials; Bernoulli polynomials; Asymptotic
expansions; Zero attractor.


%
%
\section{Introduction} \label{sec.intro}

The Appell polynomials \cite{Appell} $p_n$, $n=0,1,2,\dots$, associated 
with an entire function $g$, satisfying the condition $g(0)\neq0$, are defined 
as
\begin{equation}\label{Appell}
\sum_{n=0}^{\infty}p_n(x)\frac{z^n}{n!}  \;:=\;\frac{e^{zx}}{g(z)} \,,
\end{equation}
or, equivalently,
\begin{equation}\label{Appell_2}
p_n(x)\;=\; n![z^n]\frac{e^{xz}}{g(z)} \;=\;
                   \frac{n!}{2\pi i}\oint_\gamma \,
   \frac{e^{xz}}{g(z)}\, \frac{\mathrm{d}z}{z^{n+1}}\,,
\end{equation}
where the integration contour $\gamma$ is an oriented, index $+1$ curve
enclosing the origin and no other singularities of the
integrand.\footnote{As $g$ is required to be entire, the condition
  $g(0)\neq0$ implies the existence of a neighborhood of $z=0$ where $g$
  never vanishes.}
The symbol $[z^k]f$ denotes the coefficient of $z^k$ in the Taylor
expansion of a function $f$ about $z=0$.
These polynomials satisfy the well-known recurrence
\begin{equation}
p'_n(x) \;=\; n \, p_{n-1}(x) \,.
\end{equation}
Many important sequences of polynomials arising in analysis and 
combinatorics---including the familiar Bernoulli and Euler polynomials---
are Appell polynomials \cite{Roman}. 

For any sequence $p_n$ of polynomials, the following two questions are 
interesting:
\begin{enumerate}
\item Locate the zeros of $p_n$ (after a suitable rescaling of their 
      argument) and, in particular, determine the limiting 
      curves where they condense as $n\to +\infty$. 
\item Find the asymptotic behavior of the rescaled $p_n$ as $n\to +\infty$.
      Often there are different behaviors in different regions
      of the complex $x$-plane, separated by ``phase boundaries''.
\end{enumerate}
These two questions are closely related, as there are general theorems 
\cite{BKW1,BKW2,BKW3,sokal1} showing that the limiting curves coincide 
with the phase boundaries.
This connection arose in the pioneering work of Yang and Lee \cite{YangLee} 
on phase transitions in statistical mechanics, and was further
developed in the work of Beraha, Kahane and Weiss \cite{BKW1,BKW2,BKW3} 
and Sokal \cite{sokal1} on chromatic polynomials.  More recently,
Boyer and Goh \cite{BoyerGoh1,BoyerGoh2} have applied the above-mentioned
theorems to study the limiting curves for the Euler, Bernoulli, and Appell 
polynomials in general.  
Other recent results regarding the Appell polynomials appear in 
\cite{Buric1,Buric2}, where the authors study the asymptotics of quotients of 
gamma functions, and also in \cite{Anshelevich1,Anshelevich2,Anshelevich3}. 

The main result of the paper is to show that, the preceding integral can
be rewritten in a form that leads to interesting ways to express the
Appell polynomials and study some of their properties, in particular
the zero sets and the asymptotic expansions of the rescaled
polynomials
\begin{equation}
\pi_n(x) \;=\; p_n(nx) \,.
\label{def_rescaled_appell}
\end{equation}

We have found two representations of the Appell polynomials; each of them
has two terms: the first one is a contour integral on the complex plane,
and the second one is a sum of certain residues. Interestingly, none of
these terms is a polynomial.
The main interest of these expressions is that they show in a straightforward
way the asymptotic expansion of all of their components when $n\to+\infty$.
Another interesting feature of the the first integral term,
is that its contour is \emph{universal} for the whole family of
Appell polynomials (i.e., the same one for every acceptable choice of $g$).

This paper is organized as follows. In section~\ref{sec.appellI}, we
introduce a first integral representation for the rescaled Appell
polynomials $\pi_n(x)$ \eqref{def_rescaled_appell}.
In section~\ref{sec.other}, we give another integral
representation for these polynomials using the steepest-descent path as the
contour for the integral. In section~\ref{sec.zero_sets}, we briefly
describe how to obtain the zero attractors (in the $n\to+\infty$ limit)
for these rescaled Appell polynomials by using the asymptotic expansion of
the results of section~\ref{sec.appellI}. In section~\ref{sec.bernoulli},
we will apply the results of the previous sections to a particular (but
very important) case of the Appell polynomials: the Bernoulli polynomials
$B_n(x)$. Finally, in appendix~\ref{sec.appen}, we will collect some
interesting details about the use of the steepest-descent path as the
integration contour in section~\ref{sec.other}.

%
%
\section{An integral representation using a simple contour}
\label{sec.appellI}

\begin{figure}[htb]
\centering
\includegraphics[width=200pt]{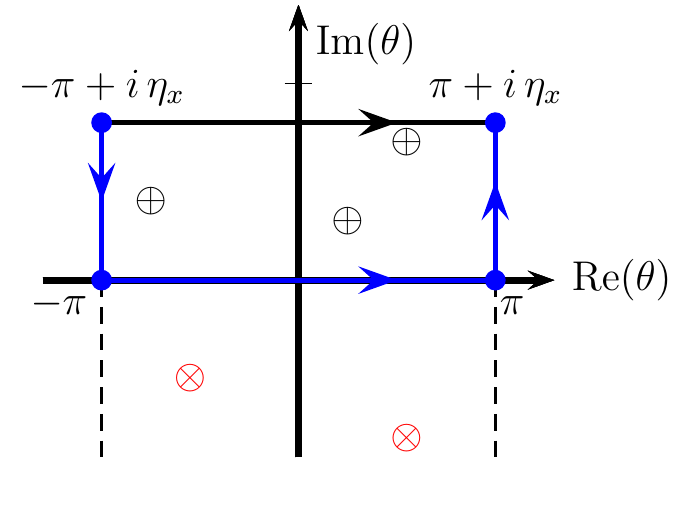}
\caption{\label{fig_appell_gen}
Contours for the integral representation \eqref{param_int} of the Appell
polynomials. The original integration contour (i.e., the
horizontal line connecting the points $-\pi+i\, \eta_x$ and $\pi+i\, \eta_x$)
can be replaced by a new one consisting of three oriented segments (two of
which will give contributions to \eqref{param_int} that cancel out due to the
periodicity of the integrand). The zeroes of $g(z)$ with positive imaginary
part (represented as $\oplus$) give rise to poles whose residues contribute
to \eqref{param_int}. The symbol $\otimes$ represents a zero of $g(z)$ with
negative imaginary part.
}
\end{figure}

Let us start by considering the integral representation \eqref{Appell_2} of
the Appell polynomials by choosing a circular integration contour $\gamma$
of radius $r$ centered at the origin of the $z$-complex plane.
For $x\neq 0$, the scaling of the integration variable $z \mapsto z/x$ gives
\begin{equation}\label{Appell_3}
\pi_n(x)\;=\; \frac{x^n n!}{2\pi i}\, \oint_{\gamma_x} \,
                 \frac{e^{nz}}{g(z/x)}\, \frac{\mathrm{d}z}{z^{n+1}}\,,
\end{equation}
where $\gamma_x$ is now a circumference of radius $r|x|$. We can adjust the
value of $r$ in such a way that no zero of $g(z/x)$ is contained within
the new integration contour. To this end we introduce $\eta_x>0$
satisfying
\begin{equation}
e^{-\eta_x} \;:=\; r|x| \;<\; r_0|x|\,,
\label{def_eta_x}
\end{equation}
where $r_0$ is the smallest modulus of the zeroes of $g(z)$.
Parametrizing the integration contour as $z=\exp(-\eta_x+i\theta)$,
we can write \eqref{Appell_3} as
\begin{equation}\label{param_int}
\pi_n(x)\;=\; \frac{x^n n!}{2\pi}\, \int_{-\pi+i\eta_x}^{\pi+i\eta_x}
\frac{\exp\big(n(e^{i\theta}-i\theta)\big)}{g(e^{i\theta}/x)}\,
\mathrm{d}\theta\,,
\end{equation}
where the integration contour is parallel to the real axis in the complex
$\theta$-plane, and goes from $-\pi+i\eta_x$ to $\pi+i\eta_x$. (See
figure~\ref{fig_appell_gen}.) Notice that we are treating $\theta$ as a
complex variable (see \cite[Chapter~6]{BenderOrszag} for a similar idea).

%
%

We can now rewrite \eqref{param_int}
as an integral along the segment of the real axis going from $-\pi$ to $\pi$
and a sum involving the residues of its integrand at the zeros $\theta_k$ of
$g(e^{i\theta}/x)$ with non-negative imaginary parts. For concreteness, we will
restrict ourselves to the case for which the zeroes with zero imaginary part
are simple. If $\theta_k$ is one of these zeroes, then the integral should
be replaced by its Principal Value (PV), and the corresponding residue
counted with a $-1/2$ coefficient. We have, hence,
\begin{equation}
\frac{\pi_n(x)}{n!} \;=\; \frac{x^n}{2\pi}\,
\left(  \mathrm{VP}\int_{-\pi}^\pi
        \frac{e^{n(e^{i\theta}-i\theta)}}{g(e^{i\theta}/x)}\mathrm{d}\theta
\;-\; 2\pi i \, \sum_{k \in \N} \Theta\left(\mathrm{Im}(\theta_k)\right) \,
\mathrm{Res}\left(
    \frac {e^{n(e^{i\theta}-i\theta)}}{g(e^{i\theta}/x)}\,;
    \theta=\theta_k\right)\right)\,,
\label{int_sum}
\end{equation}
where $\Theta \;\colon\; \R \to \R$ is the step-like function defined as:
\begin{equation}\label{def_Theta}
  \Theta(x)\;=\; \begin{cases}
    1   & \text{if $x>0$,} \\
    1/2 & \text{if $x=0$,} \\
    0   & \text{if $x<0$.}
  \end{cases}
\end{equation}
Notice that for each zero $\zeta_k$ of $g(z)$, there exists an infinite number
of zeroes of $g(e^{i\theta}/x)$ located at\footnote{
  For $z\in\mathbb{C}$ we take $\mathrm{Arg}(z)\in(-\pi,\pi]$.
}
\begin{equation}\label{zeroes}
\theta(\zeta_k,j) \;=\; \mathrm{Arg}(x\zeta_k)+2j\pi-i\log|x\zeta_k|\,,
\quad j\in\mathbb{Z}\,.
\end{equation}
Among these zeroes $\theta(\zeta_k,j)$, only a single one (namely,
$\theta_k=\theta(\zeta_k,0)$) is contained in the strip
$\mathrm{Re}(\theta)\in(-\pi,\pi]$. As $\log|x\zeta_k|\ge\log|xr_0|$, we
have $-\log|x\zeta_k|\le -\log|xr_0|<\eta_x$ and, hence, the imaginary
parts of all the singularities of the integrand in \eqref{param_int}
are strictly smaller than $\eta_x$ (in other words, they lie all below the
integration contour chosen for \eqref{param_int}). On the other hand the zeros
in \eqref{zeroes} will have positive imaginary parts if, and only if,
$|\zeta_k x|<1$,
and will lie on the integration contour (i.e., the real-$\theta$ axis with
$\theta\in(-\pi,\pi]$) if $|\zeta_k x|=1$. See figure~\ref{fig_appell_gen}
for a schematic depiction of the generic situation.

The above discussion can be summarized in the following

\begin{theorem} \label{theor.1}
The rescaled Appell polynomials $\pi_n(x)=p_n(n x)$
(cf. \eqref{def_rescaled_appell}) associated to the entire function
$g(z)$ [cf. \eqref{Appell_2}], satisfying that $g(0)\neq 0$, is equal to
\begin{equation}
\frac{\pi_n(x)}{n!} \;=\; \frac{x^n}{2\pi}\,
\left(  \mathrm{VP}\int_{-\pi}^\pi
        \frac{e^{n(e^{i\theta}-i\theta)}}{g(e^{i\theta}/x)}\, \mathrm{d}\theta
\;-\; 2\pi i \, \sum_{k \in \N} \Theta\left(\mathrm{Im}(\theta_k)\right)\,
\mathrm{Res}\left(
    \frac {e^{n(e^{i\theta}-i\theta)}}{g(e^{i\theta}/x)}\,;
    \theta=\theta_k\right)\right)\,,
\label{eq.theor.1}
\end{equation}
where $\Theta$ is given by \eqref{def_Theta}, and
$\theta_k=\theta(\zeta_k,0)$ corresponds to the zeroes $\zeta_k$ of $g(z)$
via \eqref{zeroes}.
\end{theorem}

\bigskip

\noindent
{\bf Remark}.
The expression \eqref{eq.theor.1} for the rescaled Appell polynomial $\pi_n$
is \emph{exact}. Even though the $\pi_n$ are polynomials, both terms in
\eqref{eq.theor.1} are not. This expression is very useful, as one can
obtain very easily its asymptotic expansion when $n\to+\infty$, as shown in
section~\ref{sec.zero_sets}.

%
%
\section{An integral representation using the steepest-descent path}
\label{sec.other}

\begin{figure}[htb]
\centering
\includegraphics[width=200pt]{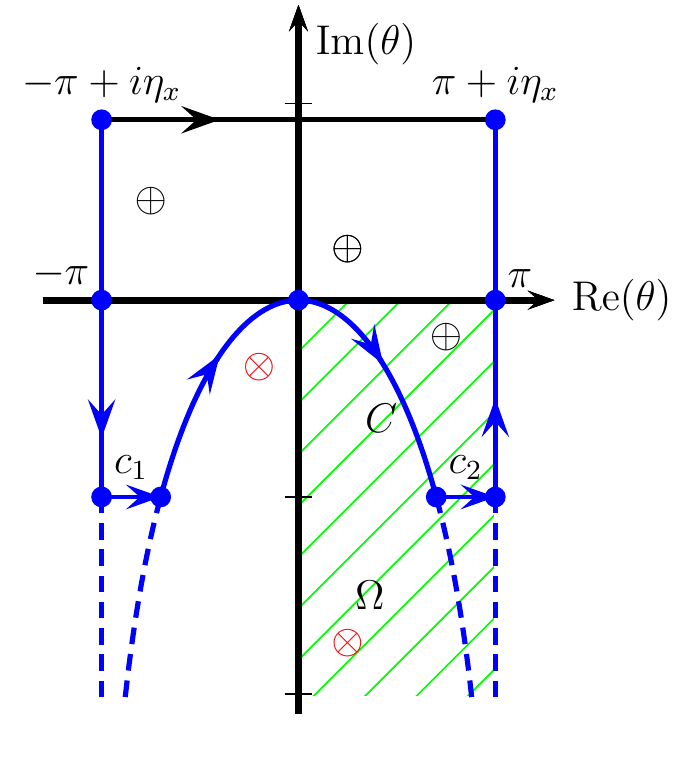}
\caption{Deformed integration contours obtained from
figure~\ref{fig_appell_gen}. The real interval $\theta \in [-\pi,\pi]$ in
figure~\ref{fig_appell_gen} has been deformed to include the
steepest-descent curve $C$, given by the locus of the points $\theta=X+ i\, Y$
satisfying eq.~\eqref{def_C} with $X\in(-\pi,\pi)$ and $Y\leq0$.
The contribution
of the vertical segments at $X=\pm \pi$ vanish due to the periodicity of the
integrand, and there are two `small' horizontal segments $c_1$ and $c_2$ that
join these vertical segments with $C$. The (negative) imaginary part of these
segments $c_{1,2}$ can be made in principle as large in absolute value as
desired. The hatched region represents the open set $\Omega$ \eqref{def_Omega}.
The notation for the zeroes of $g(z)$ is the same as in
figure~\ref{fig_appell_gen}.
}
\label{fig_appell_gen_bis}
\end{figure}

By using the representation of the Appell polynomials provided by
\eqref{param_int}, it is possible to get an alternative way to write them;
this one is very convenient to study their asymptotic behavior (in particular
for $\pi_n(x)$). This new representation is obtained by deforming the
integration contour in \eqref{param_int} to a new one given by curves
surrounding the singularities of the integrand---whose contributions can be
obtained by computing residues---, the steepest descent part $C$ of the curve
defined by the condition
\begin{equation}
\mathrm{Im}\left(e^{i\theta}-i\theta\right) \;=\; 0
\label{def_SD_curve_C}
\end{equation}
($\theta=0$ is a saddle point of $\exp(e^{i\theta}-i\theta)$), and two
straight segments $c_{1,2}$ parallel to the real axis in the $\theta$-complex
plane (see figure~\ref{fig_appell_gen_bis}). Notice that, owing to the
$2\pi$-periodicity of the integrand in the real-$\theta$ direction, we
can restrict ourselves to the strip $\mathrm{Re}(\theta)\in (-\pi,\pi]$.

If the contribution of the integration contours $c_{1,2}$ goes to zero as
the imaginary parts of their points tends to $-\infty$, and if the sum of
the residues picked up in the process of displacing the integration contour
is finite,\footnote{
    It is actually possible that the number of singularities of the
    integrand above the curve $C$ is \textit{infinite}. If that is the case,
    for the integral over $C$ to be well defined, it is necessary that the
    sum of the residues of the integrand at these points converges.
}
the integral can be extended to the full curve $C$.

The steepest descent curve\footnote{
   The imaginary axis is the \textit{steepest-ascent} curve passing through
   the origin, which is the saddle point.
}
$C$ passes through the origin of the complex $\theta$-plane. If we write
$\theta=X+iY$ with $X\in (-\pi,\pi)$, we immediately obtain a simple
implicit equation for the points of $C\setminus\{0\}$
\begin{equation}
X\, e^Y \;-\; \sin X \;=\; 0\,.
\label{def_C}
\end{equation}

%
%

A convenient way to parametrize the curve $C$ is the following:
As $C\setminus\{0\}$ is defined by the condition
$\mathrm{Im}(e^{i\theta}-i\theta)=0$, for every $\tau\neq0$ we can obtain
$\theta(\tau)$ by solving the equation
\begin{equation}\label{eq_param}
e^{i \theta(\tau)}\;-\; i\, \theta(\tau) \;=\; 1 \;-\; \tau^2\,.
\end{equation}
The $-\tau^2$ term on the r.h.s.\/  guarantees that we are indeed dealing
with the curve of steepest descent.\footnote{
   A $+\tau^2$ term would give the curve of steepest ascent.
}
Let us concentrate first on the solutions $\theta_+(\tau)$ satisfying
$\mathrm{Re}(\theta_+(\tau))\in(0,\pi)$. Notice that from each
$\theta_+(\tau)$ satisfying \eqref{eq_param}, we can obtain another
solution satisfying $\mathrm{Re}(\theta_-(\tau))\in(-\pi,0)$ as minus
its complex conjugate, i.e. $\theta_-(\tau)=-\overline{\theta_+(\tau)}$.
Notice also that $\theta(0)=0$ is also a solution of \eqref{eq_param}.

We show now that for every $\tau>0$, there is a \textit{unique} solution
$\theta_+(\tau) \in \Omega$ to \eqref{eq_param}, where the region $\Omega$
is defined by the conditions
\begin{equation}
\Omega \;=\; \left\{ \theta \in \C \;\colon\; \mathrm{Re}(\theta) \in (0,\pi)
 \wedge \mathrm{Im}(\theta) < 0 \right\} \,.
\label{def_Omega}
\end{equation}
To this end we first notice that given a solution $\theta$ of
\eqref{eq_param}, the number $z :=-e^{i\theta}$ satisfies
\begin{equation}\label{z_equation}
  z\, e^z \;=\; -e^{\tau^2-1}\,.
\end{equation}
The solutions in $z$ to this equation can be explicitly written in terms of
the different branches of the Lambert function
\cite{Corless_96,Jeffrey_96,Corless_97}.

Given one solution $\zeta$ of \eqref{z_equation}, we can immediately get a solution $\theta_\zeta$
to \eqref{eq_param} as
\begin{equation}\label{equation_z}
\theta_\zeta \;=\; i\, (1-\tau^2+\zeta)\,.
\end{equation}
For $\theta_\zeta$ to be in $\Omega$ it is necessary that
$\mathrm{Im}(\zeta)\in(-\pi,0)$. If we write $\zeta=x+iy$ with $y\in(-\pi,0)$,
equation \eqref{z_equation} is equivalent to
\begin{subequations}
\label{def_equations_z}
\begin{align}
x             &\;=\; -y\, \cot y\,, \label{equations_z_1}\\
e^{\tau^2-1}  &\;=\; -e^x\, (x\cos y \;-\; y\sin y)\,. \label{equations_z_2}
\end{align}
\end{subequations}
For pairs $(x,y)$ satisfying \eqref{equations_z_1}, the second equation
\eqref{equations_z_2} tells us that
$\tau(y)$ can be expressed as a continuous real function
$\tau \,\colon\, (-\pi,0)\to \R$ such that:
\begin{equation}\label{taudey}
\tau(y) \;=\;  \sqrt{1-y\cot y+\log\left( \frac{y}{\sin y} \right) }\,.
\end{equation}
This function is injective because
\begin{equation}
\tau'(y) \;=\; \frac{(y-\frac{1}{2}\sin2y)^2+\sin^4y}
                   {2\, y\, \tau(y)\, \sin^2y} \;<\;  0\,.
\end{equation}
By invoking Brouwer's Invariance of Domain theorem \cite[Chapter~2]{Hatcher},
we see that the map \eqref{taudey} has a continuous inverse
\begin{equation}
y \;\colon\; (0,+\infty) \to (-\pi,0) \;\colon\; y \;\mapsto\; \tau(y) \,,
\end{equation}
which is actually smooth as a consequence of the Inverse Function theorem.
As we can see, for every $\tau>0$ there is a \emph{unique} value $y(\tau)$
satisfying the equations \eqref{def_equations_z}.
From it we compute $x(\tau)=-y(\tau)\cot y(\tau)$, and find $\theta_+(\tau)$.
Notice that if $y\in (-\pi,0)$, we have $x=-y\cot y>-1$, hence, we see that
for each $\tau>0$, the \emph{unique} solution to \eqref{def_equations_z} is
contained in the open region
\begin{equation}
\mathcal{R} \;:=\; \left\{ z\in\mathbb{C} \;\colon\; \mathrm{Re}(z)>-1
                        \wedge \mathrm{Im}(z)\in(-\pi,0)\right\}\,.
\label{def_R}
\end{equation}

\begin{figure}
\centering
\begin{tabular}{cc}
\includegraphics[width=200pt]{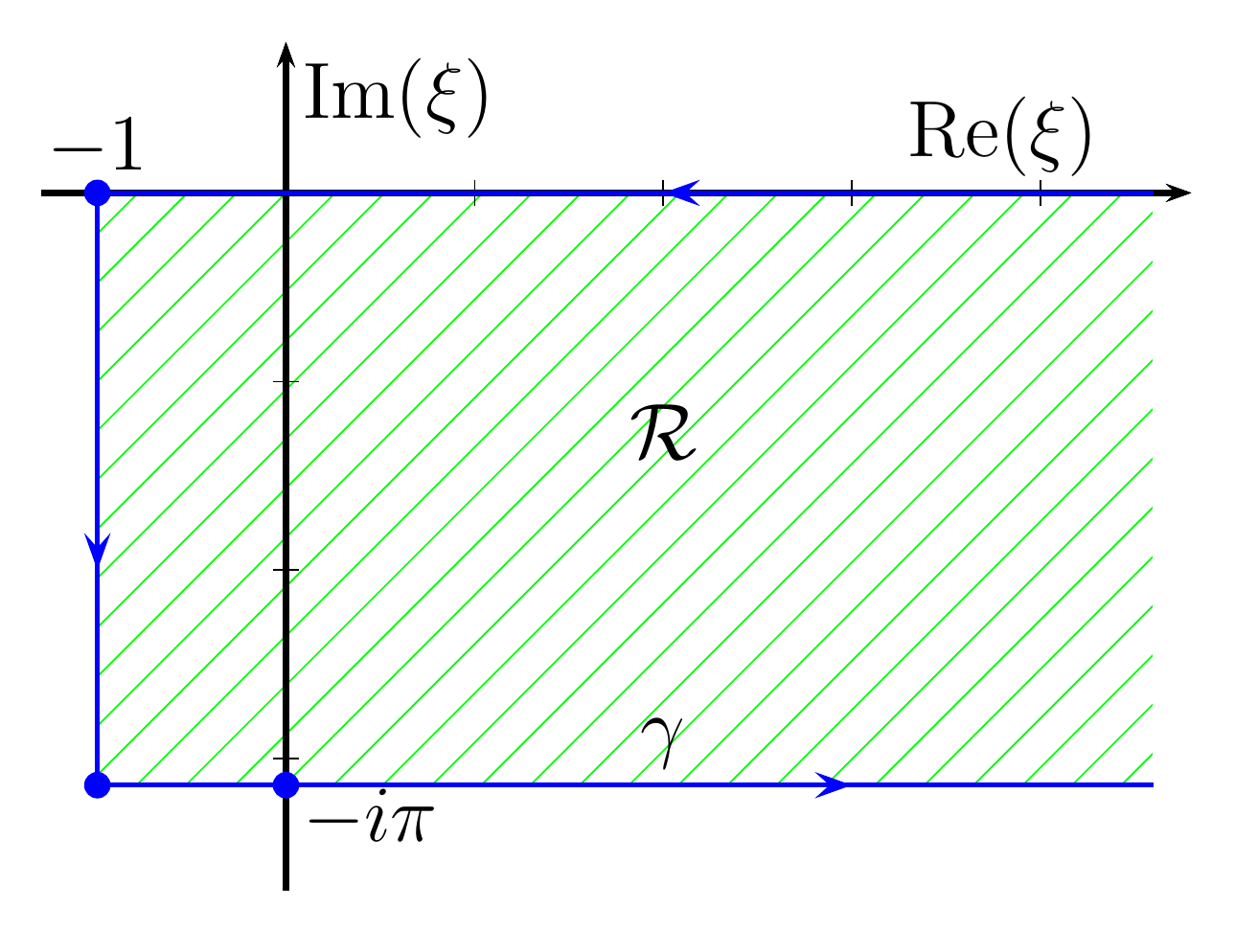} &
\includegraphics[width=170pt]{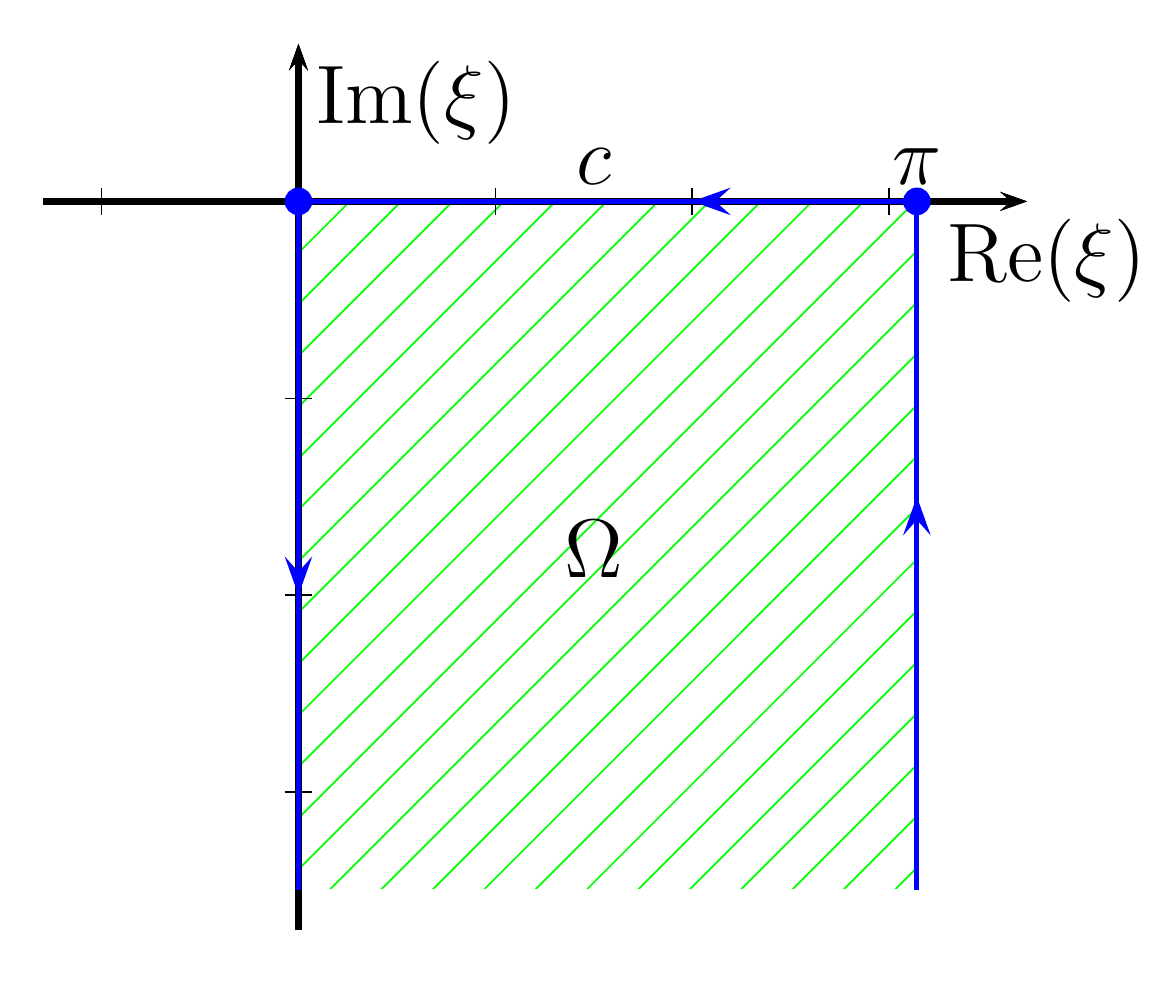} \\[2mm]
(a) & (b) \\
\end{tabular}
\caption{(a) Integration contour $\gamma$ for the integral representation of
the Lambert function that we use in the paper \eqref{Lambert1}. The contour
is the boundary of the open region $\mathcal{R}$ \eqref{def_R}.
(b) Integration contour $c$ for the integral representation of the function
$\theta_+(\tau)$ \eqref{integraltheta}. In this case, this is the boundary
of the open region $\Omega$ \eqref{def_Omega}.
(This region corresponds to the hatched region in
figure~\ref{fig_appell_gen_bis}.)}
\label{fig_contours}
\end{figure}

Now we take advantage of this fact to write a closed form expression for
the solution of \eqref{eq_param} in terms of a suitably defined branch of the
Lambert function, indeed, \eqref{equation_z} can be written as
\begin{equation}\label{sol_1}
\theta_+(\tau) \;=\; i\, \left(1-\tau^2+W\left(-e^{\tau^2-1}\right)\right)\,,
\end{equation}
where $W$ denotes the branch of the Lambert function defined by the following
integral
\begin{equation}\label{Lambert1}
  W(z) \;=\; \frac{z}{2\pi i} \, \int_\gamma \, \frac{\xi+1}{\xi e^\xi-z}\,
     \mathrm{d}\xi\,,
\end{equation}
and the integration contour $\gamma$, shown in figure~\ref{fig_contours}(a),
is the positively oriented boundary of $\mathcal{R}$ \eqref{def_R}.
The argument principle can be used to write \eqref{Lambert1} because we know
that equation \eqref{equation_z} has a \emph{unique} solution within the
region $\mathcal{R}$ delimited by $\gamma$.

%
%

It is interesting to note that, for $\tau>0$, we have
$W(-e^{\tau^2-1})=W_{-1}(-e^{\tau^2-1})$, where the $W_{-1}$ branch of the
Lambert function is defined in \cite{Corless_96,Jeffrey_96,Corless_97}.

A direct integral representation for $\theta_+(\tau)$---that can be derived
by invoking the argument principle in the domain $\Omega$ \eqref{def_Omega}--- is
\begin{equation}\label{integraltheta}
\theta_+(\tau) \;=\; \frac{1}{2\pi} \, \int_c \,
   \frac{\xi(i\xi-\tau^2)}{e^{i\xi}-i\xi-1+\tau^2}\, \mathrm{d}\xi\,,
\end{equation}
where the contour $c$ is shown in figure \ref{fig_contours}(b) and bounds the
region $\Omega$ \eqref{def_Omega}. From this integral representation, we
can immediately conclude that, given any $\tau_0>0$, there is an open
neighborhood $U_{\tau_0}\in\C$ of the half-line $\tau>\tau_0$
where $\theta_+(\tau)$ is analytic.

We can now write a parametrization for the full steepest descent curve $C$ by
defining the following function $\theta:\mathbb{R}\rightarrow\mathbb{C}$
\begin{equation}\label{def_theta}
  \theta(\tau)\;=\; \begin{cases}
         \theta_+(\tau)& \text{if $\tau>0$,} \\
         0             & \text{if $\tau=0$,} \\
         -\overline{\theta_+(-\tau)}
                      &  \text{if $\tau<0$.}
  \end{cases}
\end{equation}
It is straightforward to show that the curve $\theta(\tau)$ is smooth and
satisfies $\theta'(\tau)\neq0$ for all $\tau\in\mathbb{R}$. In fact
\begin{equation}\label{thetadot}
  \theta'(\tau) \;=\; \begin{cases}
  \displaystyle\frac{2i\tau}{i\, \theta(\tau)-\tau^2}& \text{if $\tau\neq0$,}
    \\[4mm]
  \sqrt{2}                                            &\text{if $\tau=0$.}
  \end{cases}
\end{equation}

The condition $\tau^2=1+i\, \theta-e^{i\theta}$ defines an analytic
function $\theta(\tau)$ for $\tau\in\mathbb{C}$ in a neighborhood of $\tau=0$.
To see this, notice that $(1+i\,\theta-e^{i\theta})/\theta^2$ can be
analytically extended to an entire function $h$ defined on the full
complex plane satisfying $h(0)\neq0$. This means that
\begin{equation}
\tau(\theta) \;=\; \theta\, \sqrt{\frac{1+i\, \theta-e^{i\theta}}{\theta^2}}
\end{equation}
(where we use the branch of the square root with positive real part defined
by the standard determination for the $\log$) is analytic in a neighborhood
$U_0$ of $\theta=0$, and satisfies $\tau'(0)\neq0$. By taking the $\tau_0$
that we mentioned after \eqref{integraltheta} inside the open
set $U_0$, we can actually show that the function $\theta$ \eqref{def_theta}
can be analytically extended to an open neighborhood of the real axis.

We have seen that $\tau(\theta)$ admits an analytic inverse $\theta(\tau)$
in a neighborhood of $\tau=0$. Its Taylor expansion around $\tau=0$ can be
obtained by inverting the series expansion of $\theta(\tau)$ about $\tau=0$.
However, it is interesting to mention another way to do this. The idea is
to write
\begin{equation}
\exp(\theta(\tau)) \;=\; \sum_{n=0}^\infty
      Y_n(i\theta'(0),\ldots,i\theta^{(n)}(0))\, \frac{\tau^n}{n!}\,,
\end{equation}
where $Y_n(a_1,\ldots,a_n)$ denotes the $n$-th complete Bell polynomial
(see e.g., \cite{Comtet}).
From this we get the following recurrence relation for the
$\theta^{(k)}(0)$ ($k\geq1$)
\begin{subequations}
\label{recurrence}
\begin{align}
\sqrt{2}   &\;=\; \theta'(0)\,, \label{recurrence1} \\
0          &\;=\; Y_n(i\theta'(0),\ldots,i\theta^{(n)}(0))-i\theta^{(n)}(0)
               \,. \label{recurrence2}
\end{align}
\end{subequations}
Notice that, as a consequence of the fact that the difference
$Y_n(a_1,\ldots,a_n)-a_n$ depends only on $a_1,\ldots,a_{n-1}$,
equation \eqref{recurrence2} only involves
$\theta'(0),\ldots,\theta^{(n-1)}(0)$. We give a number of terms for the
Taylor expansion of $\theta(\tau)$ around $\tau=0$ in
appendix~\ref{sec.appen}.

With the help of the parametrization of $C$ that we have introduced,
we can finally write
\begin{equation}\label{int_C}
\frac{x^n}{2\pi} \, \int_C\, \frac{\exp\big(n(e^{i\theta}-i\theta)\big)}%
   {g(e^{i\theta}/x)}\, \mathrm{d}\theta \;=\;
\frac{(ex)^n}{2\pi} \, \int_{-\infty}^\infty \,
\frac{\theta'(\tau) \, \exp(-n\tau^2) }{g(e^{i\theta(\tau)}/x)}\,
\mathrm{d}\tau\,.
\end{equation}
By using the Taylor expansions that we give in appendix~\ref{sec.appen},
it is a simple
exercise now to get as many terms of the asymptotic expansion of
\eqref{int_C} as one wants in the limit $n\rightarrow+\infty$ by relying,
for instance, on Watson's lemma (as long as $x\in \mathbb{C}$ is such that
no zeroes of the denominator lie on the real-$\tau$ axis).
Combining these expansions with the contributions of the residues picked
up in the process of deforming the integration contour to $C$,
it is straightforward to get asymptotic expansions for any Appell polynomial
evaluated at any value of the argument $x$.
The conclusion is given by

\begin{theorem} \label{theor.2}
The rescaled Appell polynomials $\pi_n(x)=p_n(n x)$
(cf. \eqref{def_rescaled_appell}) associated to the entire function
$g(z)$ [cf. \eqref{Appell_2}], satisfying that $g(0)\neq 0$, are equal to
\begin{multline}
\label{eq_theor.2}
\frac{\pi_n(x)}{n!} \;=\; \frac{x^n}{2\pi}
\left( e^n \, \mathrm{VP}\int_{-\pi}^\pi
        \frac{\theta'(\tau)\, e^{-n\tau^2}}{g(e^{i\theta(\tau)}/x)}\,
        \mathrm{d}\tau \right. \\
\left.
\;-\; 2\pi i \, \sum_{k \in \N} \Theta\left(\psi(\theta_k)\right) \,
\mathrm{Res}\left(
    \frac {e^{n(e^{i\theta}-i\theta)}}{g(e^{i\theta}/x)}\,;
    \theta=\theta_k\right)\right)\,,
\end{multline}
where $\Theta$ is given in \eqref{def_Theta}, and the function
$\psi \;\colon\; \left\{ x \in \C \;\colon\; \mathrm{Re}(x) \in (-\pi,\pi)
\right\} \to \R$ is defined by
\begin{equation}
\label{def_psi}
\psi(x) \;=\; \begin{cases}
        \mathrm{Im}(x)- \log \left( \frac{\sin \mathrm{Re}(x)}{\mathrm{Re}(x)}
        \right) & \text{if $\mathrm{Re}(x) \neq 0$, } \\
        \mathrm{Im}(x) & \text{if $\mathrm{Re}(x)=0$,}
          \end{cases}
\end{equation}
and $\theta_k=\theta(\zeta_k,0)$ correspond to the zeroes $\zeta_k$ of
$g(z)$ through \eqref{zeroes}.
\end{theorem}

\bigskip

\noindent
{\bf Remarks.} 1. The definition of the function $\psi$ \eqref{def_psi} is
motivated by the implicit equation \eqref{def_C} satisfied by the
steepest descent curve $C$. If we define $\theta= X + i Y$, the points
``above'' or on the curve $C$ are given by the expression
\begin{equation}
Y \;\ge \; \begin{cases} \log \left( \frac{\sin X}{X} \right) &
                         \text{if $X\neq 0$,} \\
                         0 & \text{if $X\neq 0$.}
           \end{cases}
\label{def_C_alternative}
\end{equation}

2. As mentioned after Theorem~\ref{theor.1}, the expression \eqref{eq_theor.2}
for the rescaled Appell polynomial $\pi_n$ is \emph{exact}, and very useful
to obtain its asymptotic expansion when $n\to +\infty$ (see next section).

%
%
\section{Zero attractors for the rescaled Appell polynomials}
\label{sec.zero_sets}

The zero attractors for some Appell polynomials have been studied by a
number of authors (see, for instance, \cite{Szego,BoyerGoh1}).
A landmark paper on this problem is \cite{BoyerGoh2}, where Boyer and Goh
provide a wealth of information about the zero sets of rescaled general
Appell polynomials (under the condition that $g$ has, at least, one zero).
The description of these sets is made easy by the application of a theorem
proved by Sokal \cite{sokal1} that narrows down the search for the points
belonging to these attractors to the determination of suitable asymptotic
approximations of their integral representations
(when $n\rightarrow+\infty$), and the study of their analyticity.
For this purpose, 
the method described in the preceding sections offers a 
quick alternative to the one used in \cite{BoyerGoh2}. 

%
%
\begin{figure}[htb]
\centering
\includegraphics[width=200pt]{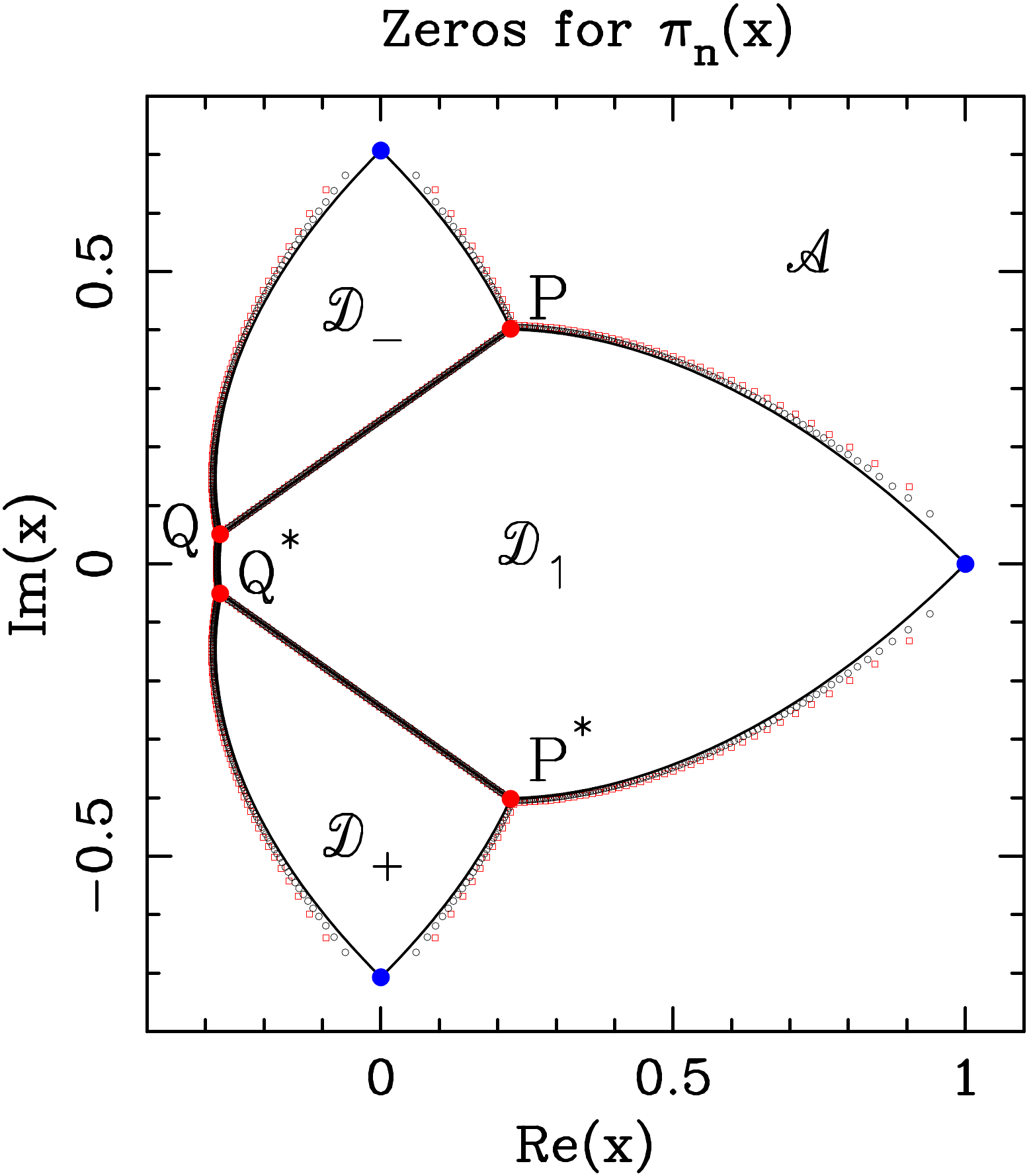}
\caption{Zero attractor for the rescaled Appell polynomials $\pi_n$ with
$g(x) =(x-1)(x^2+2)$ \cite[figure~5(a)]{BoyerGoh2}. We show the zeros for
$n=400$ (red $\Box$) and $n=1000$ (black $\circ$). The solid lines correspond
to the different branches of the zero attractor in this case. Region
$\mathcal{A}$ is where the integral \eqref{asympt_int} dominates
over the residue contributions \eqref{asympt_res}.
Region $\mathcal{D}_1$ is where the residue of the zero at $x=1$ dominates,
while regions $\mathcal{D}_\pm$ are where the residues of the zeros at
$x=\pm i \sqrt{2}$ dominate. The blue dots (at $x=1,\pm i/\sqrt{2}$) mark
the cusps corresponding to the Szeg\"o curves due to \eqref{curvas_1}, and
the straight segments $QP$ and $Q^*P^*$ are due to \eqref{curvas_2}.
The red dots mark the points where three contributions collide in
absolute values.
\label{fig_appell_example}
}
\end{figure}

In the present setting we have to determine the asymptotic behaviors of the
integral term in \eqref{int_sum} and those of each
of the residues. The exponent of the integral in \eqref{int_sum} has a
saddle point at $\theta=0$, and its asymptotic behavior
(if $g(1/x)\neq0$)\footnote{
   If $g(1/x)=0$ or, more generally, a zero of $g(e^{i\theta}/x)$ lies on
   the real axis, one should get the asymptotics of the PV of the integral
   in \eqref{asympt_int}.}
is
\begin{equation}\label{asympt_int}
\int_{-\pi}^\pi \,
\frac{\exp\big(n(e^{i\theta}-i\theta)\big)}{g(e^{i\theta}/x)} \,
\mathrm{d}\theta \;=\; \frac{e^n}{g(1/x)}\sqrt{\frac{2\pi}{n}} \,
\left(1+O\left(\frac{1}{n}\right)\right)\,,\quad
\text{as $n\rightarrow+\infty.$}
\end{equation}

The residues in \eqref{int_sum} have the form
\begin{equation}\label{residues}
\mathrm{Res}
\left(\frac{\exp\big(n(e^{i\theta}-i\theta)\big)}{g(e^{i\theta}/x)}\,;
\theta=\theta_k\right) \;=\; P_k(n)\,
\exp\big(n(e^{i\theta_k}-i\theta_k)\big)
\;=\; P_k(n)\left(\frac{e^{\zeta_k x}}{\zeta_k x}\right)^n\,,
\end{equation}
where the $P_k(n)$ are polynomials in $n$ with $x$-dependent coefficients.
The asymptotic behavior of these residues when $n\rightarrow+\infty$ is
\begin{equation}\label{asympt_res}
\mathrm{Res}\left(
\frac{\exp\big(n(e^{i\theta}-i\theta)\big)}{g(e^{i\theta}/x)}\,;
\theta=\theta_k\right) \;=\; c_k(x)\, (i\zeta_kx-i)^{p-1}\, n^{p-1}\,
\left(\frac{e^{\zeta_k x}}{\zeta_kx}\right)^n\,
\left(1+O\left(\frac{1}{n}\right)\right)\,,
\end{equation}
where $p$ denotes the order of the zero $\theta_k$, and
\begin{equation}
c_k(x)\;=\;\frac{1}{(p-1)!}\, \lim_{\theta\rightarrow\theta_k}
           \frac{(\theta-\theta_k)^p}{g(e^{i\theta}/x)}\,.
\end{equation}
For each $x\in \mathbb{C} \setminus\{0\}$, the asymptotic behavior of
$\pi_n(x)$ is then determined by the function
\begin{equation}\label{fi}
\Phi(x) \;=\; \max_{k \colon |\zeta_kx|<1}\left\{ e,
\left|\frac{e^{\zeta_kx}}{\zeta_kx}\right|\right\}\,.
\end{equation}
This is the result given in \cite{BoyerGoh2}. In order to obtain the actual
sets where the zeroes accumulate, it is necessary to find the curves
defined by conditions of the type
\begin{subequations}
\label{curvas12}
\begin{align}
\left|\frac{e^{\zeta_k x}}{\zeta_k x}\right| &\;=\; e \,, \quad \text{or,}
\label{curvas_1}
\\[2mm]
\left|\frac{e^{\zeta_k x}}{\zeta_k x}\right| &\;=\;
\left|\frac{e^{\zeta_j x}}{\zeta_j x}\right| \quad
\text{for some indexes $k\neq j$.} \label{curvas_2}
\end{align}
\end{subequations}
In the first case \eqref{curvas_1}, we obtain rotated and rescaled Szeg\"o
curves, and in the second \eqref{curvas_2}, straight lines
(see \cite{BoyerGoh2}). Notice that in the particular case of taking
$g(z)=z-\zeta$ for some $\zeta\in \mathbb{C} \setminus \{0\}$, the Appell
polynomials are directly related to the truncated exponential, and their
zero sets are rotated and dilated Szeg\"o curves (as dictated by the value
of $\zeta$).

In figure~\ref{fig_appell_example}, we show the zeros of the Appell polynomials
associated
to the choice $g(x) =(x-1)(x^2+2)$ \cite[figure~5(a)]{BoyerGoh2}, as well as
the corresponding zero attractors in the limit $n\to+\infty$.

%
%
\section{Rescaled Bernoulli polynomials}\label{sec.bernoulli}

The Bernoulli polynomials $B_n$ are the Appell polynomials for
\begin{equation}
g(z) \;=\; \frac{ e^z -1 }{z} \,,
\label{def_g_Bernoulli}
\end{equation}
so, the integral representation for the rescaled Bernoulli polynomials
$\beta_n(x) := B_n(n x)$ is given by \eqref{int_sum}/\eqref{def_g_Bernoulli}.
Notice that, as required, $g$ is an entire function with $g(0)=1 \neq 0$.

\begin{figure}[htb]
\centering
\includegraphics[width=200pt]{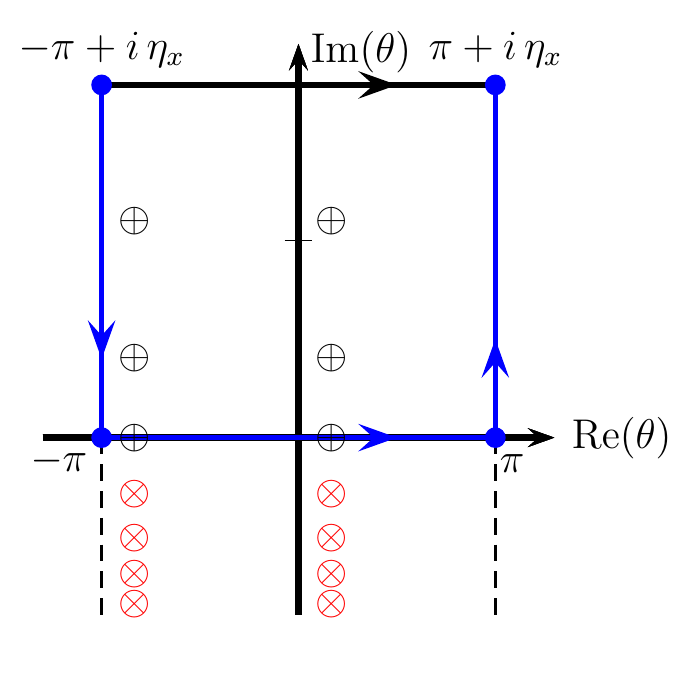}
\caption{\label{fig_bernoulli}
Generic situation for the integral representation of the rescaled
Bernoulli polynomials $\beta(x)$ \eqref{param_int}/\eqref{def_g_Bernoulli}.
This picture should be compared to figure~\ref{fig_appell_gen}. We have
chosen $x = e^{2\pi i/3}/(6\pi)$ (with $|x| < 1/(2\pi)$), and
$r=\pi < r_0=2\pi$, so that $\eta_x = -\log(r|x|) = \log(6)$ and
$k_\text{max}=3$. Notice that in this case, $\theta_3^\pm$ take real values.
The zeroes of $g$ with $k\le 7$ contributing (resp.\/ not contributing)
to the sum  in \eqref{eq.betan_final} are depicted as $\oplus$
(resp.\/ $\otimes$). They are located at
$\theta_k^+ = -5\pi/6 - i \log(k/3)$ and
$\theta_k^- =   \pi/6 - i \log(k/3)$.
}
\end{figure}

The zeros $\zeta_k$ of $g(z)$ are given by
$\zeta_k = 2\pi i k$ with $k\in \Z \setminus \{0\}$. This means that
$r_0 = 2\pi$ in this case. As $g$ is evaluated at $e^{i\theta}/x$,
the zeros of $g$ correspond to the following values of $\theta$:
\begin{equation}
\theta_k \;=\; \mathrm{Arg}(k x i) -i \, \log |2\pi k x|\,, \quad
      \text{for $k \in \Z \setminus \{0\}$}\,.
\label{def_theta_k1}
\end{equation}
In the following we prefer to work with
\begin{equation}
\theta_k^\pm  \;:=\; \mathrm{Arg}(\pm x i) -i \, \log (2\pi k |x|)\,, \quad
      \text{for $k \in \N$.}
\label{def_theta_kpm}
\end{equation}
Notice that all zeroes $\theta_k^\pm$ are simple, and
that the condition $\mathrm{Im}(\theta_k^\pm) \ge 0$ implies that
\begin{equation}
       2\pi k |x| \;\le\; 1  \,.
\end{equation}
This condition (that depends on $|x|$) will give the values of $k$ associated
to the contributing zeroes. If for some $k_\circ\in \N$,
$2\pi k_\circ |x| = 1$, the corresponding zeroes $\theta_{k_\circ}^\pm$ will
lie on the integration contour.
The maximum value of $k\in\N$ contributing to the sum in
\eqref{int_sum} is given by
\begin{equation}
k_\text{max}(x) \;=\; \left\lfloor \frac{1}{2\pi |x|} \right\rfloor \,.
\label{def.kmax}
\end{equation}
The general situation has been depicted on figure~\ref{fig_bernoulli}.

%
%

The residues of the integrand at $\theta_k^\pm =
\mathrm{Arg}(\pm x i) -i \, \log (2\pi k |x|)$, $k\in\N$
[cf. Eq.~\eqref{def_theta_kpm}] are
\begin{equation}
\mathrm{Res}\left(
\frac{\exp\big(n(e^{i\theta}-i\theta)\big)}{g(e^{i\theta}/x)};
\theta=\theta^\pm_k\right) \;=\; -\frac{i}{(2\pi k x)^n} \,
e^{\mp i n (\pi/2 -2 \pi k x)} \,.
\end{equation}
The residues for $\theta_k^+$ and $\theta_k^-$ can be combined to
give a cosine. We then get the

\begin{corollary} \label{coro.1}
The rescaled Bernoulli polynomials $\beta_n(x)=B_n(n x)$ are equal to
\begin{equation}
\frac{\beta_n(x)}{n!}\;=\; \frac{x^n}{2\pi}\, \mathrm{PV}\!\int_{-\pi}^\pi
    \frac{e^{n(e^{i\theta}-i\theta)}}{g(e^{i\theta}/x)} \,
    \mathrm{d}\theta
  -\frac{2}{(2\pi)^n}\, \sum\limits_{k\in \mathbb{N}}
  \Theta\left(\frac{1}{2\pi |x|}-k\right)\, \frac{1}{k^n} \,
  \cos\left(n\left( \frac{\pi}{2}-2\pi k x \right) \right) \,,
\label{eq.betan_final}
\end{equation}
where $g$ is given by \eqref{def_g_Bernoulli}, and $\Theta$,
by \eqref{def_Theta}.
\end{corollary}

\bigskip

\noindent
{\bf Remarks.} 1. Again, the expression \eqref{eq.betan_final}
for the rescaled Bernoulli polynomial $\pi_n$ is \emph{exact} for any $n\in\N$.

2. For any $x\in \C \setminus \{0\}$, the sum in \eqref{eq.betan_final}
is always finite. The number of terms in that sum is given by
$k_{\mathrm{max}}$ \eqref{def.kmax}.

3. This representation is useful to obtain the asymptotic expansion of
$\beta_n(x)$ as $n\to +\infty$ because the behavior of the cosine terms
is obvious, and the asymptotics of the integral in this limit (at least,
the leading term) is easy to obtain (see the remaining of this section
for a general approach).

4. It is instructive to compare this expression with the Fourier series for
the Bernoulli polynomials $B_n(x)$ \cite{Temme}. As we can see, if we rescale
back the variable $x \rightarrow x/n$, the upper limit $k_{\mathrm{max}}$
in the sum of \eqref{eq.betan_final} grows to $+\infty$ when
$n\rightarrow+\infty$. This, together with the fact that in this limit the
contribution of the integral to \eqref{eq.betan_final} is subdominant
with respect to the one of the trigonometric sum, gives as a result the
known Fourier representation for the Bernoulli polynomials.

\bigskip

\begin{figure}[htb]
\centering
\includegraphics[width=200pt]{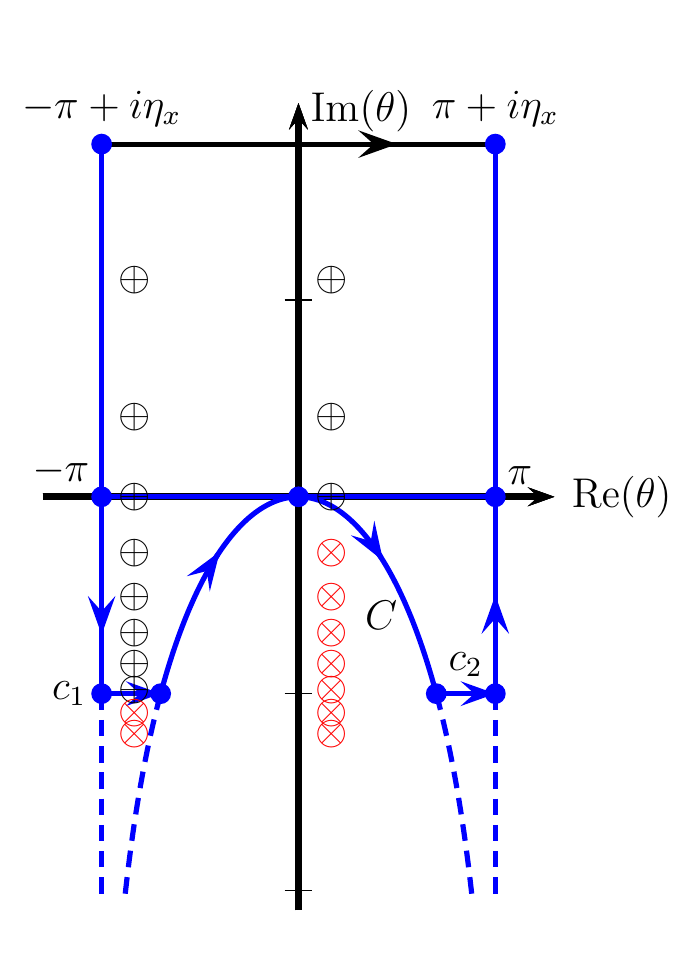}
\caption{\label{fig_bernoulli_bis}
Deformed integration contours obtained from figure~\ref{fig_bernoulli}
(in the same way as figure~\ref{fig_appell_gen_bis} was obtained from
figure~\ref{fig_appell_gen}). The new contour contains the
steepest descent curve $C$ \eqref{def_C}, as well as two small horizontal
segments $c_{1,2}$. The notation is as in
figures~\ref{fig_appell_gen_bis} and~\ref{fig_bernoulli}.
The parameters $x$, $r$, and $\eta_x$ are chosen as in
figure~\ref{fig_bernoulli}. We have shown the first $10$ zeroes
$\theta_k^\pm$. The residues that contribute to \eqref{sol_Bernoulli_1}
correspond to $k\le 8$ for $\theta^-_k$, and to $k\le 3$ for
$\theta^+_k$. In this particular example, none of these zeroes falls
on the integration contour.
}
\end{figure}

We discuss now the representation of the $\beta_n(x)$ obtained by deforming
the integration contour in \eqref{param_int} to $C$ and taking into account
the residue contributions picked up in the process (as we did in
section~\ref{sec.other}). See figure~\ref{fig_bernoulli_bis} for a depiction
of this situation.

%
%

If $x$ is not purely imaginary (i.e., if $\mathrm{Re}(x)\neq0$), a
straightforward computation gives
\begin{multline}
\frac{\beta_n(x)}{n!} \;=\;
\frac{(ex)^n}{2\pi}\,\mathrm{VP}\int_{-\infty}^{+\infty}\,
\frac{e^{-n\tau^2}\theta'(\tau)}{g(e^{i\theta(\tau)}/x)}\mathrm{d}\tau
-\frac{1}{(2\pi)^n}\, \sum_{k\in\N}
\Theta\left(\frac{\mathrm{Arg}(ix)}{2\pi\mathrm{Im}(ix)}-k\right)\,
   \frac{e^{2k\pi i n x}}{(ik)^n} \\
  - \frac{1}{(2\pi)^n}\, \sum_{k\in\mathbb{N}}
\Theta\left(\frac{\mathrm{Arg}(-ix)}{2\pi\mathrm{Im}(-ix)}-k\right)\,
\frac{e^{-2k\pi i n x}}{(-ik)^n} \,, \qquad\qquad
\label{sol_Bernoulli_1}
\end{multline}
where $g$ is given by \eqref{def_g_Bernoulli} and $\Theta$, by
\eqref{def_Theta}. Notice that the upper limits in the sums of the r.h.s.\/
of \eqref{sol_Bernoulli_1} are finite in this case ($\mathrm{Re}(x)\neq0$),
but may be different from each other for generic choices of
$x\in\mathbb{C}\setminus \{0\}$. (See figure~\ref{fig_bernoulli_bis}.)

\bigskip

\noindent
{\bf Remark.} The argument of the function $\Theta$ in \eqref{sol_Bernoulli_1}
comes directly from the generic condition \eqref{def_C_alternative} for a
zero $\theta_k$ to be ``above'' or on the steepest descent curve $C$
\eqref{def_C}, and the particular form for these zeros for the Bernoulli case
\eqref{def_theta_kpm}.

\bigskip

For $x=iq$ with $q>0$, we find
\begin{multline}
\frac{ \beta_n(x)}{n!} \;=\;
\frac{(ex)^{n}}{2\pi}\,\mathrm{VP}\int_{-\infty}^{+\infty}
\frac{e^{-n\tau^2}\theta'(\tau)}{g(e^{i\theta(\tau)}/x)}\, \mathrm{d}\tau
- \frac{1}{(2\pi i)^n} \, \mathrm{Li}_n\big(e^{2\pi i n x}\big) \\
-\frac{1}{(2\pi)^n}\, \sum_{k\in\mathbb{N}}
\Theta\left(\frac{1}{2\pi \mathrm{Im}(x)}-k\right)\,
\frac{e^{-2\pi i k n x}}{(-ik)^n} \,,
\label{sol_Bernoulli_2}
\end{multline}
where $\mathrm{Li}_n$ denotes the polylogarithm, and $\Theta$ is given by
\eqref{def_Theta}. The expression for
$x=iq$ for $q<0$ is obtained from the previous one by complex conjugation.

\medskip

\noindent
Several comments are in order now

\bigskip

\noindent
1) Whenever the upper limits in the sums on the r.h.s.\/ of
\eqref{sol_Bernoulli_1} are smaller than one, the $\beta_n(x)/n!$ are
given \textit{exactly} by integral on the r.h.s.\/ of \eqref{sol_Bernoulli_1},
which is written in a form that makes it specially easy to study its
asymptotic behavior when $n\rightarrow+\infty$.

The first two terms of the asymptotic expansion of the integral that
appears on the r.h.s.\/ of \eqref{sol_Bernoulli_1} in the limit
$n\rightarrow+\infty$ can be obtained by using the expressions given in
appendix~\ref{sec.appen}. In the case $\exp(1/x)-1\neq0$, we have
\begin{multline}
\frac{(ex)^{n}}{2\pi}\,
\mathrm{VP}\int_{-\infty}^{+\infty}
\frac{e^{-n\tau^2}\theta'(\tau)}{g(e^{i\theta(\tau)}/x)}\, \mathrm{d}\tau
\;=\;
\frac{e^n x^{n-1}}{(e^{1/x}-1)\, \sqrt{2\pi n}} \\
\times \left(1-
\frac{x^2+e^{2/x}(6-12x+x^2)-2e^{1/x}(x^2-6x-3)}{12x^2(e^{1/x}-1)^2}\,
\frac{1}{n}+O\left(\frac{1}{n^2}\right)\right)\,,
\label{asympt_int_1}
\end{multline}
whereas, for $1/x=2\pi i m$ with $m\in\mathbb{Z}\setminus\{0\}$,
the asymptotic expansion when $n\to+\infty$ of the integral on the r.h.s.\/ of
\eqref{sol_Bernoulli_2} is
\begin{subequations}
\begin{align}
\frac{(ex)^{n}}{2\pi}\,&\mathrm{VP}\int_{-\infty}^{+\infty}
\frac{e^{-n\tau^2}\theta'(\tau)}{g(e^{i\theta(\tau)}/x)}\, \mathrm{d}\tau
\nonumber \\
& \;=\;
\frac{(ex)^n}{2\sqrt{2\pi n}}\,
\left(\frac{2}{3}-\frac{1}{x}+
\left(\frac{1}{270}+\frac{1}{12x}-\frac{1}{6x^2}\right)\,\frac{1}{n}
+O\left(\frac{1}{n^2}\right)\right) \\[2mm]
&\;=\;\left(\frac{e}{2\pi i m}\right)^n \, \frac{1}{2\sqrt{2\pi n}}\,
   \left(\frac{2}{3}-2\pi i m+
\left(\frac{1}{270}+\frac{m \pi i}{6}+\frac{2m^2\pi^2}{3}\right)\, \frac{1}{n}
+O\left(\frac{1}{n^2}\right)\right)\,.
\end{align}
\end{subequations}

\medskip

\noindent
2) The contributions of the curves $c_{1,2}$ (see
figure~\ref{fig_appell_gen_bis}) can be easily seen to go to zero as they
are displaced in the direction of the negative imaginary axis. This is a
consequence of the rapid fall off of the term
$\exp(n(e^{i\theta}-i\theta))$. In the case $\mathrm{Re}(x)\neq0$, only
finite sums are involved, so the convergence of the integral in
equation \eqref{sol_Bernoulli_1} is guaranteed. If $\mathrm{Re}(x)=0$,
some care is needed, as we have to deal with an infinite number of
contributions from singularities located at points with real parts equal
to $\pm \pi$. In this case, the sum will be proportional to the
polylogarithm $\mathrm{Li}_n(e^{\pm 2\pi i n x})$, which is finite in general.
The best way to proceed is to restrict oneself to contours
$c_{1,2}$ with imaginary part equal to $-\log\big((2k+1)\pi |x|\big)$
and bound the integrals on those specific contours.

\medskip

\noindent
3) As far as the asymptotics of the $\beta_n(x)$ in the limit
$n\rightarrow+\infty$ is concerned, the previous representations tell us
the relative contributions of the different terms. Notice, for instance,
that the denominators in the sums on the r.h.s.\/ of \eqref{sol_Bernoulli_1}
exponentially suppress the contributions of the terms with $k>1$ (if present).
The asymptotic behavior of the integral term in \eqref{sol_Bernoulli_1}
is controlled by an exponential factor $(ex)^n$. If this factor turns out to
outweigh the contributions of the $k=1$ terms of the these sums, all
the terms in the asymptotic expansion of the integral will dominate over
the ones coming from the sums.

%
%

\section*{Acknowledgments}

The authors wish to thank Alan Sokal for correspondence and useful comments. 
This work has been supported by the Spanish MINECO research grant
FIS2014-57387-C3-3-P and the
FEDER/Ministerio de Ciencia, Innovaci\'on y Universidades-Agencia Estatal
de Investigaci\'on/FIS2017-84440-C2-2-P grant.

\appendix

\section[Parametrization of the steepest descent curve C]%
        {Parametrization of the steepest descent curve $\bm{C}$}
\label{sec.appen}

By solving the recurrence relation \eqref{recurrence}, we can easily get
as many terms of the Taylor expansion of $\theta(\tau)$ around $\tau=0$
as we need. The terms up to $O(\tau^9)$ are
\begin{equation}\label{theta_Taylor}
\theta(\tau)\;=\;\sqrt{2} \tau -\frac{i \tau^2}{3} -\frac{\tau^3}{9 \sqrt{2}}+
\frac{2 i \tau^4}{135}+ \frac{\tau^5}{540 \sqrt{2}}+ \frac{4 i \tau^6}{8505}+
\frac{139 \tau^7}{340200 \sqrt{2}} -\frac{2 i \tau^8}{25515}+O(\tau^{9})\,.
\end{equation}

For an entire function $g$ such that $g(0)\neq0$ and $g(1/x)\neq0$, we get
to $O(\tau^3)$
\begin{multline}
\frac{\theta'(\tau)}{g\big(\exp(i\theta(\tau))/x\big)} \;=\;
\frac{\sqrt{2}}{g_0(x)}-\frac{2(xg_0(x)+3g_1(x))}{3xg_0^2(x)}\tau \\
  -\frac{x^2g_0^2(x)+12g_1^2(x)-6g_0(x)g_2(x)}{3\sqrt{2}x^2g_0^3(x)}\tau^2+
  O(\tau^3)\,, \qquad
\label{desarrollo_g}
\end{multline}
where we use the shorthand notation $\displaystyle g_n(x):=g^{(n)}(1/x)$.

If $g(1/x)\neq0$ and $x\neq0$ is such that no singularities of the integrand
lie on the integration contour $C$, the first terms of the asymptotic
expansion for \eqref{int_C} in the limit $n\rightarrow+\infty$ are
\begin{multline}
\frac{x^n}{2\pi}\, \int_C
\frac{\exp\big(n(e^{i\theta}-i\theta)\big)}{g(e^{i\theta}/x)}\,
\mathrm{d}\theta
\;=\;\frac{(ex)^n}{\sqrt{2\pi n}}\,
\left(\frac{1}{g_0(x)}  \right. \\
- \left.\left(\frac{1}{12g_0(x)}+\frac{g_1^2(x)}{x^2g_0^3(x)}-
\frac{g_2(x)}{2x^2g_0^2(x)}\right)\,\frac{1}{n}
+O\left(\frac{1}{n^2}\right)\right)\,.
\label{asympt_C}
\end{multline}
If $x$ is such that some of the singularities (simple poles) of the
integrand lie on the integration contour $C$ away from $\theta=0$,
a straightforward argument leads to the conclusion that the preceding
expression \eqref{asympt_C} remains valid for the PV  of the integral
appearing on its l.h.s.

Finally, if $x$ is such that $\theta=0$ is a simple pole of the integrand,
we have
\begin{multline}
\frac{x^n}{2\pi}\,\mathrm{PV}\int_C
\frac{\exp\big(n(e^{i\theta}-i\theta)\big)}{g(e^{i\theta}/x)}\,
\mathrm{d}\theta \;=\;
\frac{(ex)^n}{2\sqrt{2\pi n}}\; \left(-\frac{4xg_1(x)+3g_2(x)}{3g_1^2(x)}
\right. \\
  +\frac{92x^3g_1^3(x)+135g_2^3(x)-180g_1(x)g_2(x)g_3(x)+
  45g_1^2(x)\big(x^2g_2(x)+g_4(x)\big)}{540x^2g_1^4(x)}\,\frac{1}{n} \\
+ \left.O\left(\frac{1}{n^2}\right)\right)\,. \qquad\qquad\qquad\qquad
\label{asympt_Cpole}
\end{multline}

\end{document}